\documentclass{article}
\usepackage{amsmath, amsthm}
\usepackage{graphicx}

\newtheoremstyle{theorem}
  {10pt}		  
  {10pt}  
  {\sl}  
  {\parindent}     
  {\bf}  
  {. }    
  { }    
  {}     
\theoremstyle{theorem}
\newtheorem{theorem}{Theorem}

\newtheoremstyle{defi}
  {10pt}		  
  {10pt}  
  {\rm}  
  {\parindent}     
  {\bf}  
  {. }    
  { }    
  {}     
\theoremstyle{defi}
\newtheorem{definition}[theorem]{Definition}

\begin{document}

\title{Discrete differential geometry of tetrahedrons and encoding of local protein structure}

\author{Naoto Morikawa\\
Genocript\\
e-mail: nmorika@genocript.com}

\maketitle

\begin{abstract}
Local protein structure analysis is informative to protein structure analysis and has been used successfully in protein structure prediction and others. 
Proteins have recurring structural features, such as helix caps and beta turns, which often have strong amino acid sequence preferences. 
And the challenges for local structure analysis have been identification and assignment of such common short structural motifs. 

This paper proposes a new mathematical framework that can be applied to analysis of the local structure of proteins, 
where local conformations of protein backbones are described using differential geometry of folded tetrahedron sequences. 
Using the framework, we could capture the recurring structural features without any structural templates, 
which makes local structure analysis not only simpler, but also more objective. 
Programs and examples are available from \texttt{http://www.genocript.com}. 
\\


\end{abstract}


\section{Introduction}
Protein is a sequence of amino acids, which folds into a unique three-dimensional structure in nature. And one could identify proteins with polygonal chains obtained by connecting the center of adjacent amino acids. Since the functional properties of proteins are largely determined by the structure, protein structure analysis is crucial to the study of proteins. 

Local protein structure analysis is informative to protein structure analysis and has been used successfully in protein structure prediction and others. 
Proteins have recurring structural features, such as helix caps and beta turns, which often have strong amino acid sequence preferences. 
And the challenges for local structure analysis have been identification and assignment of such common short structural motifs (\cite{BB}, \cite{BEH}, \cite{RRW}, \cite{SSL}, \cite{US}). Identification involves description of protein backbone conformation and definitions of the structural motifs. And assignment is not a trivial task, due to the variations observed in nature when compared to ideal ones.

In this paper, we introduce a new differential geometrical approach for local structure analysis. 
As for differential geometrical description, a lot of works on the surface of protein molecules are known (to name a few, \cite{BER}, \cite{CCL}). But protein backbone structure is usually studied via classification (\cite{TA}, \cite{RF}) and differential geometrical approach has been rarely taken so far.

\begin{figure}[tb]
\includegraphics{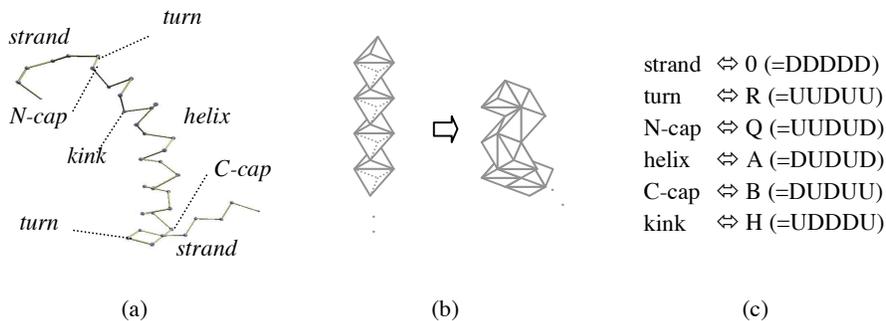}
\caption{Introduction. (a): Protein (transferase 1RKL) and its structural features. 
(b): Folding of a tetrahedron sequence. (c): Local features and the corresponding 5-tile codes (see section \ref{sec:encoding}).}
\label{fig1}
\end{figure}

One of the few is the early work of \cite{RS} which described protein backbones as polygonal chain, where each line segment corresponds to the virtual-bond between consecutive $\alpha$-carbons. In contrast, we describe local conformation of protein backbones using folded tetrahedron sequences (Figure \ref{fig1} (b)). 

As for the shape of protein backbones, \cite{EM} proposed the notion of alpha-shape and \cite{DLO} examined geometric restrictions on polygonal protein chains. Moreover, \cite{T} reviewed topological knots in protein structure.

\section{Differential geometry of triangles}
For simplicity, we first consider the differential geometry of triangles. 

\subsection{Basic ideas}
\begin{figure}[tb]
\includegraphics{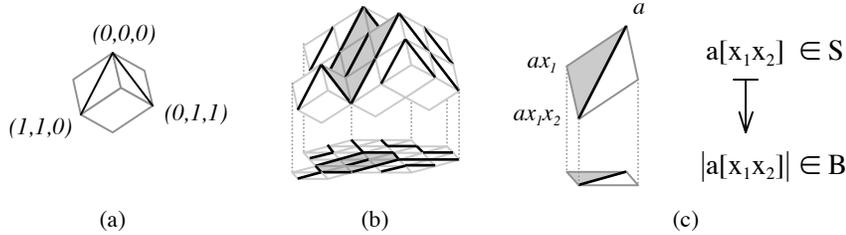}
\caption{Basic ideas. (a): Division of facets of unit cube $[0,1]^3$. (b): ``Peaks and valleys'' of cubes and its projection on a hyperplane.
(c): Projection $\pi$ from the collection $S$ of all the ``slant'' triangles to the collection $B$ of all the ``flat'' triangles.}
\label{fig2}
\end{figure}

Let's consider unit cube $[0,1]^3$ in the three-dimensional Euclidean space $\mathbf{R}^3$ and divide each of the three facets which contain $(0,0,0)$ into two triangles along diagonal, as shown in figure \ref{fig2} (a). 
Then, if we pile the cubes up in the direction of $(-1,-1,-1)$, we would obtain ``peaks and valleys'' of cubes, where the division of the facets of each cube makes up a division of the surface of the peaks and valleys (figure \ref{fig2} (b) top). 
And a ``flow'' of triangles in $\mathbf{R}^2$ is obtained by projecting the surface onto a hyperplane, (figure \ref{fig2} (b) bottom). For example, the grey  ``slant'' triangles on the surface specify the closed trajectory of the grey ``flat'' triangles on the hyperplane. 

In the following, we use monomial notation to denote points and triangles in $\mathbf{R}^3$. 
That is, we denote point $(l, m, n) \in \mathbf{R}^3$ by monomial $x_1^lx_2^mx_3^n \in \mathbf{Z}[x_1, x_2, x_3]$. 
And the triangle of vertices $(l, m, n)$, $(l+1, m, n)$, $(l+1, m+1, n)$ $\in \mathbf{R}^3$ are denoted by $x_1^lx_2^mx_3^n[x_1x_2]$. 
For example, $a[x_ix_j]$ is the triangle of vertices $a$, $ax_i$, and $ax_ix_j \in \mathbf{Z}[x_1, x_2, x_3]$ (figure \ref{fig2} (c)).

\subsection{Tangent bundle over flat triangles}

\begin{figure}[tb]
\includegraphics{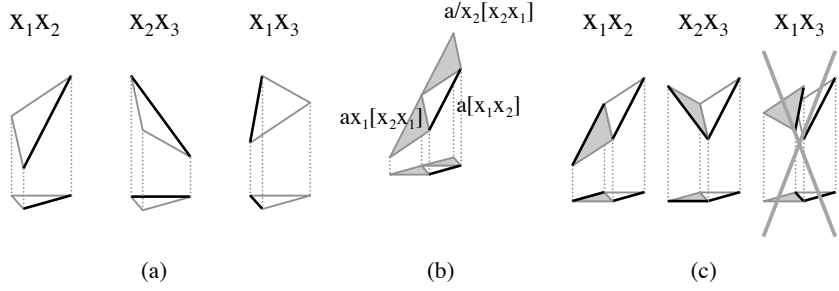}
\caption{Differential structure. (a): Gradient values of slant triangles over $|a[x_1x_2]| \in B$. 
From left to right, $a[x_1x_2]$ whose gradient is $x_1x_2$, $ax_1[x_2x_3|$ whose gradient is  $x_2x_3$, and $a/x_3[x_3x_1]$ whose gradient is $x_1x_3$.
(b): The local trajectory specified by $a[x_1x_2] \in S$ at $|a[x_1x_2]| \in B$, $|ax_1[x_2x_1]|$ (downward) and $|a/x_2[x_2x_1]|$ (upward) $\in B$. 
(c): Smoothness condition at $|a[x_1x_2]| \in B$ (white) specified by $a[x_1x_2] \in S$. The next triangle $|ax_1[x_2x_1]| \in B$ (grey) could assume either $x_1x_2$ or $x_2x_3$ as gradient.}
\label{fig3}
\end{figure}

Let $\pi$ be the projection of the collection $S$ of all the slant triangles onto the collection $B$ of all the flat triangles along direction $(-1, -1, -1)$, where the image of $a[x_ix_j] \in S$ is denoted by $|a[x_ix_j]|$ (figure \ref{fig2} (c)). 
Then, projection $\pi$ induces tangent bundle-like structure $TB$ over $B$, where the gradient of slant triangles are defined as follows:

\begin{definition}
The \textit{gradient} $Da[x_ix_j]$ of $a[x_ix_j] \in S$ is monomial $x_ix_j \in \mathbf{Z}[x_1, x_2, x_3]$. 
In particular, there is a one-to-one correspondence between $TB$ and $\{x_1x_2, x_1x_3, x_2x_3\} \times B$. 
And we indicate the gradient value over a flat triangle by a bold edge as shown in figure \ref{fig3} (a). 
\end{definition}

For example, slant triangles $a[x_1x_2]$, $ax_1[x_2x_3]$, and $a/x_3[x_3x_1] \in S$ are projected onto the same flat triangle $|a[x_1x_2]| \in B$ 
and their gradients are $x_1x_2$, $x_2x_3$, and $x_1x_3$ respectively (figure \ref{fig3} (a)).

Then, a gradient value over a flat triangle specifies a local trajectory at the flat triangle as follows:

\begin{definition}
The \textit{local trajectory} defined by $a[x_ix_j] \in S$ at $|a[x_ix_j]| \in B$ is the three consecutive flat triangles 
$\{  |ax_i[x_jx_i]|,  |a[x_ix_j]|, |a/x_j[x_jx_i]|  \} \subset B$. 
As figure \ref{fig3} (b) shows, these are the adjacent triangles connected along the direction of the bold edge of $|a[x_ix_j]|$.
And the local trajectory is specified uniquely by the gradient of $a[x_ix_j]$.
\end{definition}

Now we impose a kind of ``smoothness condition'' as shown in figure \ref{fig3} (c). 
That is, each flat triangle assume one of two gradient values, which are determined naturally by the gradient of the preceding triangle. 
Suppose that the gradient at current triangle $|a[x_1x_2]| \in B$ is $x_1x_2$ and the gradient at next triangle $|ax_1[x_2x_1]| \in B$ is $x_1x_3$. Then, two flat triangles $|a[x_1x_2]|$ and $|ax_1[x_2x_1]|$ are separated by the bold edge of $|ax_1[x_2x_1]|$ (figure \ref{fig3} (c) right).  In this case, we permit either $x_1x_2$ or $x_2x_3$ as gradient of $|ax_1[x_2x_1]|$.

\begin{figure}[tb]
\includegraphics{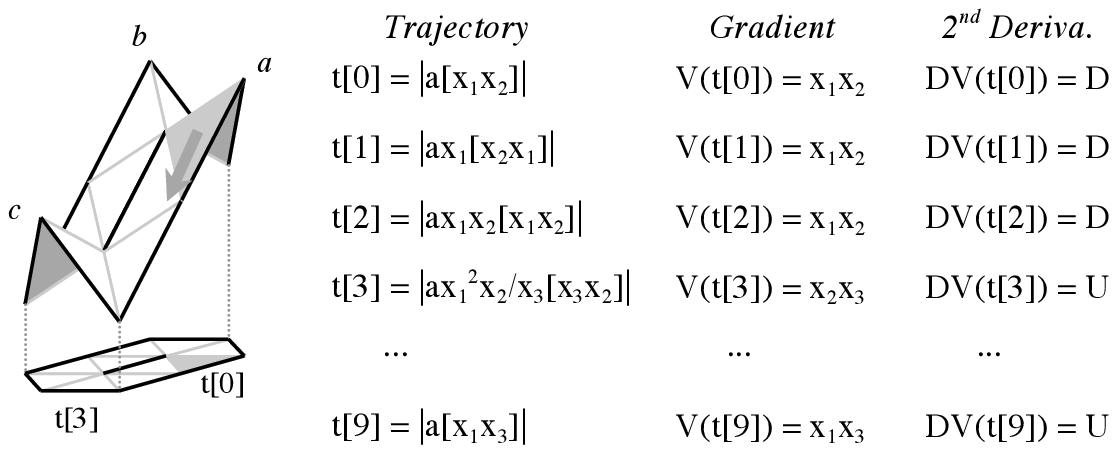}
\caption{Closed trajectory of the vector field specified by three peaks $a$, $b=a/x_2$, and $c=ax_1^2x_2/x_3$ $\in \mathbf{Z}[x_1, x_2, x_3]$.}
\label{fig4}
\end{figure}

As an example, let's consider the peaks and valleys shown in figure \ref{fig2} (b), which is specified by three peaks 
$a$, $b=a/x_2$, and $c=ax_1^2x_2/x_3$ $\in \mathbf{Z}[x_1, x_2, x_3]$ (figure \ref{fig4}). 
Peaks and valleys define a ``smooth'' vector field $V$ on $B$ by the following mapping: 
\[
V: B \rightarrow \{ x_1x_2, x_1x_3, x_2x_3 \}, \ V(|a[x_ix_j]|) := x_ix_j,
\]
where $a[x_ix_j] \in S$ is the slant triangle on the surface of the peaks and valleys over $|a[x_ix_j]| \in B$.

Let's start from triangle $|a[x_1x_2]|$ (grey) and move downward: $t[0]=|a[x_1x_2]|$ and $V(t[0])=x_1x_2$. 
Then, gradient $V(t[0])$ specifies local trajectory $\{  |ax_1[x_2x_1]|$, $|a[x_1x_2]|,$ $|a[x_1x_3]|  \}$ at $t[0]$ (Note that $|a/x_2[x_2x_1]| = |a[x_1x_3]|$). 
Since we move downward, next triangle $t[1]$ is $|ax_1[x_2x_1]|$ and we obtain $V(t[1])=x_1x_2$. 
Then, gradient $V(t[1])$ specifies local trajectory $\{ |ax_1x_2[x_1x_2]|,$ $|ax_1[x_2x_1]|,$ $|a[x_1x_2]|\}$ at $t[1]$. And next triangle $t[2]$ is $|ax_1x_2[x_1x_2]|$. Continuing the process, we obtain the closed trajectory of length $10$. 

\subsection{Encoding of the shape of trajectories}

Finally, let's consider variation of gradient along a trajectory. 
Thanks for the smoothness condition, variation of gradient, i.e., the ``second derivative'', along a triangle trajectory is given as binary valued sequence.

\begin{definition}
The \textit{derivative} $DV$ of vector field $V$ along trajectory $\{t[i]\}$ is defined as follows:
\[
DV: B \rightarrow \{U, D\}, \  DV(t[i]) := \left\{
\begin{array}{l l}
DV(t[i-1]) \  & \text{if $V(t[i]) = V(t[i-1])$} \\
 - DV(t[i-1]) \  & \text{otherwise},
\end{array}
\right.
\]
where $-U:=D$ and $-D:=U$. In words, \textit{change value if the gradient changes}.
\end{definition}

As an example, let's consider the trajectory of figure \ref{fig4} again. 
First, set any initial value: $DV(t[0]) = D$. 
Then, since the first two triangles $t[0]$ and $t[1]$ have the same gradient, $DV(t[1])$ is also $D$. 
The value of the second derivative is $D$ until $t[3]$, where it changed to $U$ since the gradient of $t[2]$ is different from that of $t[3]$. 

Continuing the process, we obtain a binary sequence of length $10$, $DDDU$$D$
$U$$UUDU$, which describes the shape of the trajectory.

\begin{figure}[tb]
\includegraphics{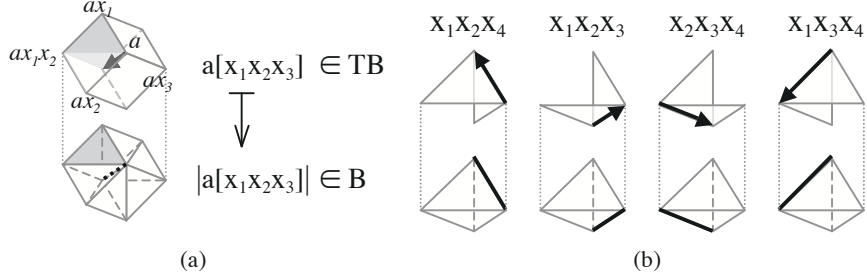}
\caption{Tangent bundle (a): Projection $\pi$ from the collection $S$ of all the slant tetrahedrons to the collection $B$ of all the flat tetrahedrons.
(b): Gradient values of slant tetrahedrons over $|a[x_1x_2x_3]| \in B$. 
From left to right, $a/x_4[x_4x_1x_2]$ whose gradient is $x_1x_2x_4$, $a[x_1x_2x_3|$ whose gradient is  $x_1x_2x_3$, $ax_1[x_2x_3x_4]$ whose gradient is $x_2x_3x_4$, and $ax_1x_2[x_3x_4x_1]$ whose gradient is $x_1x_3x_4$. 
The arrows of slant tetrahedrons indicate the direction of ``down'' in $\mathbf{R}^4$}
\label{fig5}
\end{figure}

\section{Differential geometry of tetrahedrons}

Similarly we obtain a flow of tetrahedrons in $\mathbf{R}^3$ by considering peaks and valleys of $4$-cubes in $\mathbf{R}^4$. In this case, each trajectory of tetrahedrons could be obtained by folding a tetrahedron sequence which satisfies the following conditions (figure \ref{fig1} (b)) : 
    (i) Each tetrahedron consists of four short edges and two long edges,
        where the ratio of the length is $\sqrt{3}/2$ and
    (ii) Successive tetrahedrons are connected via a long edge and have
        a rotational freedom around the edge.
In particular, we could compute the differential structure on a trajectory without considering $4$-cubes.

\subsection{Tangent bundle over flat tetrahedrons}

Let's consider $4$-cube $[0,1]^4$ in the four-dimensional Euclidean space $\mathbf{R}^4$. 
Then, the facets of $4$-cubes are three-dimensional unit cubes and we divide each of the four facets which contain $(0,0,0,0)$ into six tetrahedrons along diagonal, as shown in figure \ref{fig5} (a) top. 

In the following, we denote point $(k, l, m, n) \in \mathbf{R}^4$ by monomial $x_1^kx_2^lx_3^mx_4^n \in \mathbf{Z}[x_1, x_2, x_3,x _4]$. 
And the tetrahedron of vertices $(k, l, m, n)$, $(k+1, l, m, n)$, $(k+1, l+1, m, n)$, $(k+1, l+1, m+1, n)$ $\in \mathbf{R}^4$ are denoted by $x_1^kx_2^lx_3^mx_4^n[x_1x_2x_3]$. For example, $a[x_ix_jx_k]$ is the tetrahedron of vertices $a$, $ax_i$, $ax_ix_j$, and $ax_ix_jx_k \in \mathbf{Z}[x_1, x_2, x_3,x _4]$. 

Let $\pi$ be the projection of the collection $S$ of all the slant tetrahedrons onto the collection $B$ of all the flat tetrahedrons along direction $(-1, -1, -1, -1)$, where the image of $a[x_ix_jx_k] \in S$ is denoted by $|a[x_ix_jx_k]|$ (figure \ref{fig5} (a)). 
Then, projection $\pi$ induces tangent bundle-like structure $TB$ over $B$, where the gradient of slant tetrahedrons are defined as follows:

\begin{definition}
The \textit{gradient} $Da[x_ix_jx_k]$ of $a[x_ix_jx_k] \in S$ is monomial $x_ix_jx_k \in \mathbf{Z}[x_1, x_2, x_3, x_4]$. 
In particular, there is a one-to-one correspondence between $TB$ and $\{x_1x_2x_3, x_1x_2x_4, x_1x_3x_4, x_2x_3x_4\} \times B$. 
And we indicate the gradient value over a flat triangle by a bold edge as shown in figure \ref{fig5} (b), where arrows of slant tetrahedrons indicate the direction of ``down'' in $\mathbf{R}^4$. 
\end{definition}

For example, slant tetrahedrons $a/x_4[x_4x_1x_2]$, $a[x_1x_2x_3|$, $ax_1[x_2x_3x_4]$, and $ax_1x_2[x_3x_4x_1] \in S$ are projected onto the same flat tetrahedron $|a[x_1x_2x_3]|  \in B$ and their gradients are $x_1x_2x_4$, $x_1x_2x_3$, $x_2x_3x_4$, and $x_1x_3x_4$ respectively (figure \ref{fig5} (b)).

Then, a gradient value over a flat tetrahedron specifies a local trajectory at the flat tetrahedron as follows:

\begin{figure}[tb]
\includegraphics{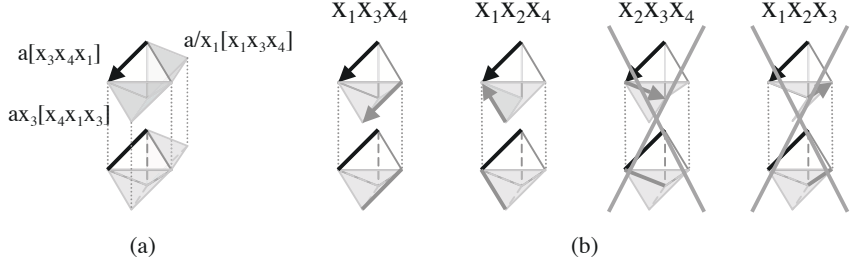}
\caption{Local trajectory. (a): The local trajectory specified by $a[x_3x_4x_1] \in S$ at $|a[x_3x_4x_1]| \in B$, $|ax_3[x_4x_1x_3]|$ (downward) and $|a/x_1[x_1x_3x_4]|$ (upward) $\in B$. 
(b): Smoothness condition at $|a[x_3x_4x_1]| \in B$ (white) specified by $a[x_3x_4x_1] \in S$. The next tetrahedron $|ax_3[x_4x_1x_3]| \in B$ (grey) could assume either $x_1x_3x_4$ or $x_1x_2x_4$ as gradient.}
\label{fig6}
\end{figure}

\begin{definition}
The \textit{local trajectory} defined by $a[x_ix_jx_k] \in S$ at $|a[x_ix_jx_k]| \in B$ is the three consecutive flat tetrahedrons 
$\{  |ax_i[x_jx_kx_i]|, |a[x_ix_jx_k]|,$ \\
$|a/x_k[x_kx_ix_j]|  \}$ $\subset B$. 
As figure \ref{fig6} (a) shows, these are the adjacent tetrahedrons connected along the direction of the bold edge of $|a[x_ix_jx_k]|$.
And the local trajectory is specified uniquely by the gradient of $a[x_ix_jx_k]$.
\end{definition}

Now we impose a kind of ``smoothness condition'' as shown in figure \ref{fig6} (b). 
That is, each flat tetrahedron assume one of two gradient values, which are determined naturally by the gradient of the preceding tetrahedron. 
Suppose that the gradient at current tetrahedron $|a[x_3x_4x_1]| \in B$ is $x_1x_3x_4$ and the gradient at next tetrahedron $|ax_3[x_4x_1x_3]| \in B$ is either $x_2x_3x_4$ or $x_1x_2x_3$. Then, the bold edges of the two flat tetrahedrons $|a[x_3x_4x_1]|$ and $|ax_3[x_4x_1x_3]|$ are not connected smoothly as shown in figure \ref{fig6} (b). In this case, we permit either $x_1x_3x_4$ or $x_1x_2x_4$ as gradient of $|ax_3[x_4x_1x_3]|$.

\begin{figure}[tb]
\includegraphics{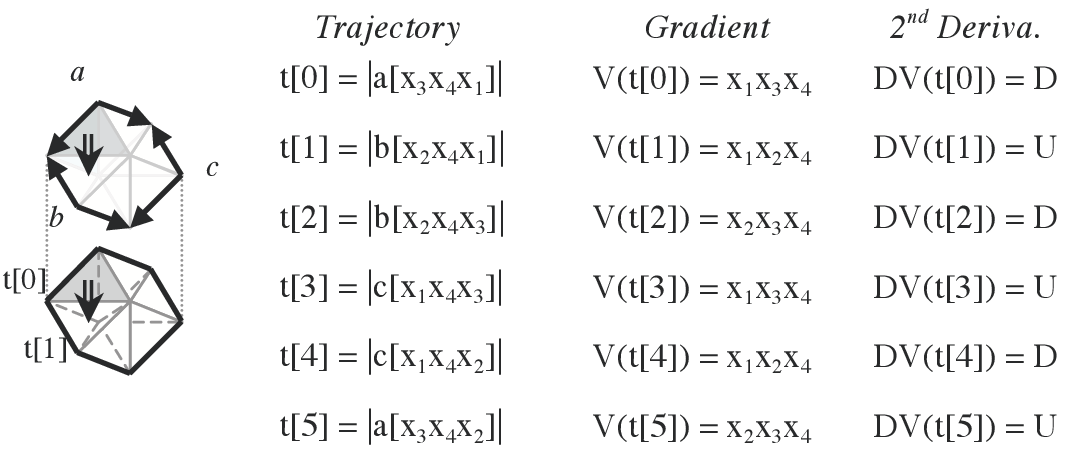}
\caption{A closed trajectory of the vector field specified by three peaks $a=x_1x_2x_4$, $b=x_1x_3x_4$, and $a=x_2x_3x_4$ $\in \mathbf{Z}[x_1, x_2, x_3, x_4]$.}
\label{fig7}
\end{figure}

As an example, let's consider a closed trajectory of peaks and valleys specified by three peaks $a=x_1x_2x_4$, $b=x_1x_3x_4$, and $a=x_2x_3x_4$ $\in \mathbf{Z}[x_1, x_2, x_3, x_4]$ (figure \ref{fig7}). 
Peaks and valleys define a ``smooth'' vector field $V$ on $B$ by the following mapping: 
\[
V: B \rightarrow \{ x_1x_2x_3, x_1x_2x_4, x_1x_3x_4, x_2x_3x_4 \}, \ V(|a[x_ix_jx_k]|) := x_ix_jx_k,
\]
where $a[x_ix_jx_k] \in S$ is the slant tetrahedron on the surface of the peaks and valleys over $|a[x_ix_jx_k]| \in B$.

Let's start from tetrahedron $|a[x_3x_4x_1]|$ (grey) and move downward: $t[0]=|a[x_3x_4x_1]|$ and $V(t[0])=x_1x_3x_4$. 
Then, gradient $V(t[0])$ specifies local trajectory $\{  |b[x_2x_4x_1]|,$ $|a[x_3x_4x_1]|,$ $|a[x_3x_4x_2]|  \}$ at $t[0]$. 
Since we move downward, next tetrahedron $t[1]$ is $|b[x_2x_4x_1]|$ and we obtain $V(t[1])=x_1x_2x_4$. 
Then, gradient $V(t[1])$ specifies local trajectory $\{ |b[x_2x_4x_3]| $ $|b[x_2x_4x_1]|$, $|a[x_3x_4x_1]| \}$ at $t[1]$. 
And next tetrahedron $t[2]$ is $|b[x_2x_4x_3]|$. Continuing the process, we obtain the closed trajectory of tetrahedrons. 

Note that the trajectory could be obtained by folding the tetrahedron sequence mentioned above (figure \ref{fig1} (b)).

\subsection{Encoding of the shape of trajectories}

Thanks for the smoothness condition again, variation of gradient, i.e., the ``second derivative'', along a tetrahedron trajectory is also given as binary valued sequence.

\begin{definition}
The \textit{derivative} $DV$ of vector field $V$ along trajectory $\{t[i]\}$ is defined as follows:
\[
DV: B \rightarrow \{U, D\}, \  DV(t[i]) := \left\{
\begin{array}{l l}
DV(t[i-1]) \  & \text{if $V(t[i]) = V(t[i-1])$} \\
 - DV(t[i-1]) \  & \text{otherwise},
\end{array}
\right.
\]
where $-U:=D$ and $-D:=U$. In words, \textit{change value if the gradient changes}.
\end{definition}

As an example, let's consider the trajectory of figure \ref{fig7} again. 
First, set any initial value: $DV(t[0]) = D$. 
Then, since the first two tetrahedrons $t[0]$ and $t[1]$ have different gradient values, $DV(t[1])$ is $U$. 
The third tetrahedron $t[2]$ assumes yet another gradient and value of the second derivative changes to $U$. 
Continuing the process, we obtain a binary sequence of length six, $DU$$DU$$DU$, which describes the shape of the trajectory.

\section{Encoding of local protein structure}\label{sec:encoding}

Now let's encode local protein structure using variation of gradient along a trajectory of tetrahedrons. 

\begin{figure}[tb]
\includegraphics{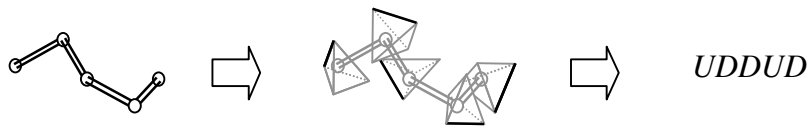}
\caption{Encoding of local structure. 
Left: A polygonal chain which represents the structure of an amino acid fragment to be encoded. 
Middle: Approximation by a folded tetrahedron sequence.
Right: Variation of gradient along the folded tetrahedron sequence.}
\label{fig8}
\end{figure}

To study the local structure of a protein, i.e., polygonal chain, we consider all the amino acid fragments of length \textit{five} occurred in the protein. (It will turn out that length five is enough to detect local features.) 
And polygonal chains are approximated by folded tetrahedron sequences to detect their local features, 
where we permit translation and rotation during the folding process to absorb irregularity of the structure (figure \ref{fig8}). 

Each fragment is approximated by a folded tetrahedron sequence of length five, starting from the middle point amino acid, say A. 
And variation of gradient along the sequence is computed to encode its structural features. 
We call the resulting \{$D$, $U$\}-valued sequence of length five the \textit{$5$-tile code} of A. 

\subsection{Encoding algorithm}

In the following, we will explain the algorithm of ``tetrahedron folding with translation and rotation.'' 
As an example, let's consider the polygonal chain $AA[-2]$-$AA[-1]$-$AA[0]$-$AA[1]$-$AA[2]$ of figure \ref{fig9} (a) and compute the 5-tile code of $AA[0]$ using a sequence of five tetrahedrons $T[-2]$-$T[-1]$-$T[0]$-$T[1]$-$T[2]$.

\subsubsection{Step 1}

Align tetrahedron $T[0]$ (white) with amino acid $AA[0]$ and set initial values (figure \ref{fig9} (b)). 
In this example, the gradient and second derivative of $T[0]$ is $x_1x_2x_4$ and $D$ respectively.

Then, the initial positions of adjacent tetrahedrons $T[\pm 1]$ (grey) are also determined, which are moved to the positions of $AA[\pm 1]$ respectively later.

\subsubsection{Step 2}

\begin{figure}[tb]
\includegraphics{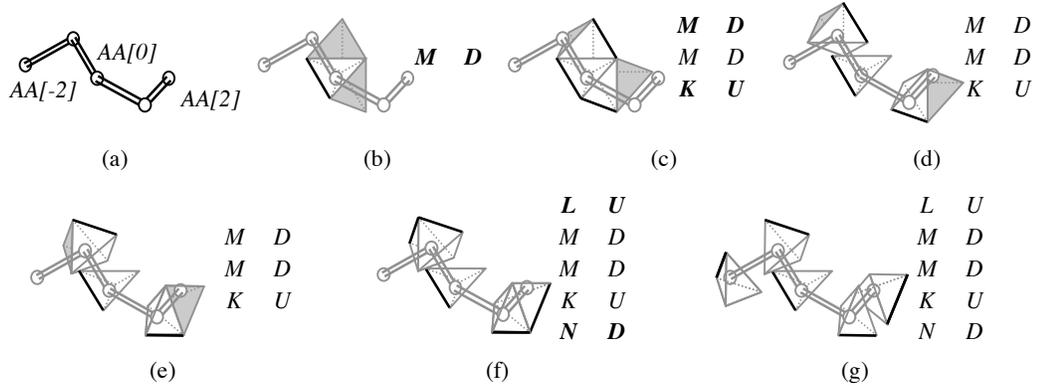}
\caption{Algorithm of the 5-tile coding. 
(a): Polygonal chain $AA[-2]$-$AA[-1]$-$\cdots$-$AA[2]$ which represents the structure of an amino acid fragment to be encoded. 
(b): Step 1. (c): Step 2. (d): Step 3. (e): Step 4. (f): Step 5. (g) Step 6. 
The character strings show the corresponding sequences of gradients (left) and second derivatives (right), 
where top are the those of $T[-2]$ and bottom are those of $T[2]$. $K$, $L$, $M$, and $N$ stand for $x_2x_3x_4$, $x_1x_2x_3$, $x_1x_2x_4$, and $x_1x_3x_4$ respectively.}
\label{fig9}
\end{figure}

Assign gradient to adjacent tetrahedrons $T[\pm 1]$ considering the direction of $AA[\pm 2]$ respectively (Figure \ref{fig9} (c)). 
For example, tetrahedron $T[1]$ could assume $x_1x_2x_4$ or $x_2x_3x_4$ as its gradient. And the next tetrahedron (grey) becomes closer to $AA[2]$ if $x_2x_3x_4$ is assumed.
Thus, the gradient of $T[1]$ is $x_2x_3x_4$ and its second derivative is $U$ since the gradients of $T[0]$ and $T[1]$ are different. 
In the same way,  $T[-1]$ is assigned $x_1x_2x_4$ and $D$ as its gradient and second derivative respectively.

Note that the initial positions of adjacent tetrahedrons $T[\pm 2]$ (grey) are also determined, which are moved to the positions of $AA[\pm 2]$ respectively later.

\subsubsection{Step3}

Translate tetrahedrons $T[\pm 1]$ to the positions of $AA[\pm 1]$ respectively (Figure \ref{fig9} (d)). 
Adjacent tetrahedrons $T[\pm 2]$ (grey) are also moved with $T[\pm 1]$ respectively. 

\subsubsection{Step4}

Rotate tetrahedrons $T[\pm 1]$ at the positions of $AA[\pm 1]$ so that the bold edges become parallel to the direction from $AA[0]$ to $AA[\pm 2]$ respectively (Figure \ref{fig9} (e)). 
Adjacent tetrahedrons $T[\pm 2]$ (grey) are also moved with $T[\pm 1]$ respectively.

\subsubsection{Step5}

Assign gradient to adjacent tetrahedrons $T[\pm 2]$ considering the direction of $AA[\pm 2]$ respectively (figure \ref{fig9} (f)). 
For example, tetrahedron $T[2]$ could assume $x_1x_3x_4$ or $x_2x_3x_4$ as its gradient. And the next tetrahedron (not shown) becomes closer to $AA[2]$ if $x_1x_3x_4$ is assumed.
Thus, the gradient of $T[2]$ is $x_1x_3x_4$ and its second derivative is $D$ since the gradients of $T[1]$ and $T[2]$ are different. 
In the same way,  $T[-2]$ is assigned $x_1x_2x_3$ and $U$ as its gradient and second derivative respectively.

\subsubsection{Step6}

Translate tetrahedrons $T[\pm 2]$ to the positions of $AA[\pm 2]$ respectively (figure \ref{fig9} (g)). 
And we have obtained binary sequence $UDDUD$, the 5-tile code of $A[0]$, which describes the shape of the amino acid fragment shown in figure \ref{fig9} (a).

\subsection{One-letter representation of 5-tile codes}

To save space, we use numerals and alphabets to denote 5-tile code $C_1C_2C_3C_4C_5$. 
First, compute the value $Y$ of the code which is defined as follows: 
$Y = 2^4*C_1' + 2^3*C_2' + 2^2*C_3' + 2*C_4' + C_5'$, 
where $C_i'=1$ if $C_i$ is equal to $U$ and $C_i'=0$ if not.
Then, assign the number to the code if the value $Y$ is less than 10. Otherwise, assign the $(Y-9)$-th alphabet to the code.

For example, $DDDUU$ corresponds to binary number $00011$ and $Y=3$. Thus, $3$ is assigned to the code. 
On the other hand, $DUDUD$ corresponds to binary number $01010$ and $Y=10$. Thus, the first alphabet $A$ is assigned to the code. 

\subsection{Example: transferase 1RKL}

The local structure of transferase 1RKL shown in figure \ref{fig1} (a) is encoded as follows:
\[
\begin{array}{cccc}
{\tt MISDEQLNSL} &{\tt AITFGIVMMT} &{\tt LIVIYHAVDS} & {\tt TMSPKN}\\
{\tt \cdot \cdot 000RQAAA} &{\tt AAAAHAAAAA} & {\tt AAAAAAB0R0} & {\tt 0000 \cdot \cdot},
\end{array}
\]
where the top row shows the amino acid sequence of 1RKL and the bottom shows the corresponding 5-tile codes.
As you see, we could capture recurring structural features without any structural templates (figure \ref{fig1} (c)).

In previous works, common short structural motifs (structural templates) of proteins are often identified by clustering a set of representative protein fragments, using unsupervised machine learning. Thus, identification and assignment of such motifs has been the challenges for local structure analysis. 
And, as a result,  their methods could not recognize new local structural features nor structural distortions. 

On the other hand, there is no need for identification and assignment of structural templates in our method since we don't use them at all. 
And the $5$-tile codes could detect both new local features and structural distortions because they are computed directly from atomic coordinates.

\end{document}